\documentclass[12pt]{article}
\usepackage[mathscr]{eucal}
\usepackage{amssymb,amsmath}
\font\cmc=cmcsc10  scaled \magstep2

\newcommand\vk{\vskip}
\newcommand\hk{\hskip}

\newcommand\al{\alpha}
\newcommand\be{\beta}
\newcommand\iy{\infty}

\newcommand\la{\lambda}
\newcommand\ga{\gamma}

\newcommand\varep{\varepsilon}
\newcommand\rg{\rightarrow}

\newcommand\De{\Delta}

\newcommand\ov{\overset}
\newcommand\und{\underset}
\newcommand\no{\noindent}

\newcommand\ovl{\overline}
\newcommand\for{\forall}
\newcommand\col{\colon\hk-.5em}

\newtheorem{Proof.}{\it Proof.}

\begin{document}
\vbox to .5truecm{}

\begin{center}
\cmc A Characterization of Extremal Sets in Hilbert Spaces
\end{center}
\vk.3cm
\begin{center}
by Viet NguyenKhac
\footnote{Institute of Mathematics, P.O.Box 631 Bo Ho, 
10000 Hanoi, Vietnam;\ e-mail:\ nkviet@thevinh.ncst.ac.vn}  
\ \& Khiem NguyenVan
\footnote{Department of Mathematics and Informatics, Hanoi University of Education, Cau Giay dist., Hanoi, Vietnam.}
\end{center}
\begin{center} 
Hanoi Institute of Mathematics \& Hanoi University of Education
\end{center}
\vk.1cm
\begin{center}
\end{center}

\normalsize
\begin{abstract} {We give a characterization of extremal sets in Hilbert spaces that generalizes a classical theorem of H. W. E. Jung. We investigate also the behaviour of points near to the circumsphere of such a set with respect to the Kuratowski and Hausdorff measures of non-compactness.}
\end{abstract}
\vk.5cm					
\begin{center} 
\textbf {1.\ \ Introduction}
\end{center}

Let $(X,\|\cdot\|)$ be a Banach space. For a non-empty bounded subset $A$ of $X$ and a non-empty subset $B$ of $X$ we shall use the following notations:\ \ $d(A)\col=$\ sup $\{\|x-y\|:\ \ x, y\in A\}$ -- the diameter of $A$;\ $r_B(A)\col=\und{y\in B}{\text {inf}}\ \und{x\in A}{\text {sup}}\ \|x-y\|$ 
-- the relative Chebyshev radius of $A$ with respect to $B$, in particular $r(A)\col=r_{\ovl{co} A}(A)$ with $\ovl{co} A$ denoting the closed convex hull of $A$;\ $C_B(A)\col=\{y\in B:\ \und{x\in A}{\text {sup}}\ \|x-y\|=r_B(A)\}$ -- the set of Chebyshev centers of $A$ in $B$.  
\vk.2cm
The Jung constant of $X$ is defined by $J(X)\col=$ sup $\{r_X(A):\ A\subset X,\ \ {\text {with}}\ \ d(A)=1\}$. In conection with uniform normal structure one considers also the following important geometric constant - the self-Jung constant of $X$ ({\it cf.} \cite{Byn}):\ \ $J_s(X)\col=$ sup $\{r(A):\ A\subset X,\ \ {\text {with}}\ \ d(A)=1\}$. In the case $X$ is an inner product space, it is known that $C_X(A)$ consists of a unique point which belongs to the closed convex hull $\ovl{co} A$ of $A$ ({\it cf.} \cite{Gar}). Hence $J(X)=J_s(X)$ in this case. Classical Jung's theorem states that for $X=E^n$ -- an $n$-dimensional Euclidean space we have $J(E^n)=J_s(E^n)=\sqrt{\dfrac{n}{2(n+1)}}$\ \ (\cite{Jun}, {\it cf.} \cite{Dgk}). Furthermore if $X=H$ -- a Hilbert space (the infinite-dimensional case), then $J(H)=J_s(H)=\dfrac{1}{\sqrt{2}}$\ \ (\cite{Rou}, {\it cf.} \cite{Ber}, \cite{Byn}, \cite{Das}). 
\vk.2cm
\textbf{Definition.}\quad We say that a bounded subset $A$ of $X$ consisting of at least two points is an extremal set (resp. self-extremal set), if $r_X(A)=J(X) d(A)$\ (resp. $r(A)=J_s(X) d(A)$).
\vk.2cm
Note that in the case $X=E^n$, or $H$, two these notions coincide, so in this case we shall speak simply of extremal sets. From the second part of the mentioned theorem of Jung one knowns that a bounded subset $A$ of $E^n$ is extremal if and only if $A$ contains a regular $n$-simplex with edges of length $d(A)$. In the case $X=H$ a partial result in this direction was obtained by Gulevich (\cite{Gul}) who showed that if $A$ is a relatively compact set in a Hilbert space with $d(A)>0$, then $r(A)<\dfrac{1}{\sqrt{2}} d(A)$. Therefore if $A$ is an extremal set in a Hilbert space, then $A$ is not relatively compact.
\vk.2cm
Our aim in this paper is to give a characterization of extremal sets in Hilbert spaces which is a genaralization of the second part of Jung's theorem.
\vk.5cm   
\textbf{Main Theorem.}\quad {\it Let $A$ be an extremal set in a Hilbert space $H$ with $d(A)=\sqrt{2}$. Then $\chi(A)=1$ and for every $\varep\in (0,\sqrt{2})$, for every positive integer $p$ there exists a\ $p$-simplex $\De$ with its vertices in $A$ and each edge of $\De$ has length not less than $\sqrt{2}-\varep$.

Conversly if $d(A)=\sqrt{2}$ and for every $\varep\in (0,\sqrt{2})$, for every positive integer $p$ there exists a\ $p$-simplex $\De$ with its vertices in $A$ such that the length of each edge of $\De$ is not less than $\sqrt{2}-\varep$, then $A$ is an extremal set.}
\vk.2cm
In the above formulation $\chi(A)$ denotes the Hausdorff measure of non-compactness of $A$, {\it i.e.} the infimum of positive $r$ such that $A$ can be covered by a finite number of balls with radius $r$ and with centers in $X$. Besides, based on an observation of \cite{Das} (``Mushroom Lemma") we prove also a result on the behaviour of points near to the circumsphere of such a set with respect to the measure of non-compactness which says roughly that the main contribution to the measure of non-compactness comes from that part of the extremal set.
\vk.5cm
\begin{center}
\textbf {2.\ \ Measures of non-compactness of extremal sets}
\end{center}

\textbf{Theorem\ 1.}\quad {\it Let $A$ be an extremal set in a Hilbert space $H$ with $r(A)=1$. Then we have $\al(A)=\sqrt{2}$.}
\vk.2cm
Here $\al(A)$ denotes the Kuratowski measure of non-compactness, that is the infimum of positive $d$ such that $A$ can be covered by a finite number of sets of diameter $d$.
\vk.2cm
{\it The first proof}.\quad From the assumption $r(A)=1$ it follows that for each integer number $n\ge 2$ we have $\und{x\in A}{\bigcap}\ B(x,1-\dfrac{1}{n})=\emptyset$,\ where $B(x,r)$ denotes the closed ball centered at $x$ with radius $r$ which is weakly compact since $H$ is reflexive. Hence there exist  $x_{p_{n-1}+1},\ x_{p_{n-1}+2},\ \cdots,\ x_{p_n}$ in $A$ such that $\und{i=p_{n-1}+1}{\ov{p_n}{\bigcap}}\ B(x_i,1-\dfrac{1}{n})=\emptyset$ (with convention $p_1=0$). 
\vk.2cm
Setting $A_n\col=\{x_{p_{n-1}+1},\ x_{p_{n-1}+2},\ \cdots,\ x_{p_n}\}$\ we denote the Chebyshev center of $A_n$ in $H$ by $c_n$ and let $r_n\col=r(A_n)$, then $r_n>1-\dfrac{1}{n}$.
\vk.2cm
Let $S(c,r)$ denote the sphere with center $c$ and radius $r$. From the proof of classical Jung's theorem one knowns that $A\cap S(c_n,r_n)\ne\emptyset$ and $c_n\in co(A_n\cap S(c_n,r_n))$. So there exist $y_{q_{n-1}+1},\ y_{q_{n-1}+2},\ \cdots,\ y_{q_n}$ in $A_n\cap S(c_n,r_n)$ (with convention $q_1=0$) and positive numbers $t_{q_{n-1}+1},\ t_{q_{n-1}+2},\ \cdots,\ t_{q_n}$ such that
  $$c_n=\sum_{q_{n-1}<i\le q_n}\ t_i y_i,\qquad \sum_{q_{n-1}<i\le q_n}\   t_i=1.$$

We claim that $\al(\{y_{q_{n-1}+1},\ y_{q_{n-1}+2},\ \cdots,\ y_{q_n}\}_{n=2}^\iy)=\sqrt{2}$. Suppose on the contrary $\al(\{y_{q_{n-1}+1},\ y_{q_{n-1}+2},\ \cdots,\ y_{q_n}\}_{n=2}^\iy)<\sqrt{2}$. Then one can choose $\varep_0\in (0,\sqrt{2})$ satisfying $\al(\{y_{q_{n-1}+1},\ y_{q_{n-1}+2},\ \cdots,\ y_{q_n}\}_{n=2}^\iy)\le\sqrt{2}-\varep_0$, and so subsets $D_1,\ D_2,\ \cdots,\ D_m$ of $H$ with $d(D_i)\le\sqrt{2}-\varep_0$ for every $i=1, 2, \cdots, m$ such that $\{y_{q_{n-1}+1},\ y_{q_{n-1}+2},\ \cdots,\ y_{q_n}\}_{n=2}^\iy\subset\und{i=1}{\ov{m}{\bigcup}}\ D_i$. There exists at least one set among $D_1,\ D_2,\ \cdots,\ D_m$, say $D_1$ with the property that there are infinitely many $n$ satisfying

  $$\sum_{i\in J_n}\ t_i\ge \dfrac{1}{m}\eqno{(1)}$$

\no where

  $$J_n\col=\{i\in [q_{n-1}+1,q_n]:\ \ y_i\in D_1\}.$$
\vk.2cm
For each $n$ satisfying (1) and fixed $j\in J_n$ we have

  \begin{align} & \sum_{q_{n-1}<i\le q_n}\ t_i \|y_i-y_j\|^2=\sum_{q_{n-1}<i\le q_n}\ t_i\|y_i-c_n+c_n-y_j\|^2=\notag\\
  &{\quad} =\sum_{q_{n-1}<i\le q_n}\ t_i\big( \|y_i-c_n\|^2+\|y_j-c_n\|^2-2(y_i-c_n,y_j-c_n)\big)=\notag\\
  &{\quad} =2r_n^2-2\big(\sum_{q_{n-1}<i\le q_n}\ t_i y_i-c_n,y_j-c_n\big)=\notag\\
  &{\quad} =2r_n^2>2\big(1-\dfrac{1}{n}\big)^2>2-\dfrac{4}{n},\notag
 \end{align}

\no where $(\cdot,\cdot)$ denotes the inner product in $H$.
\vk.2cm
On the other hand

  $$\begin{aligned} \sum_{q_{n-1}<i\le q_n}\ t_i \|y_i-y_j\|^2 &=\sum_{i\in J_n}\ t_i\|y_i-y_j\|^2+\sum_{q_{n-1}<i\le q_n,\ i\notin J_n}\ t_i \|y_i-y_j\|^2\le\\
  &\le (\sqrt{2}-\varep_0)^2\ \sum_{i\in J_n}\ t_i+2\big(1-\sum_{i\in J_n}\ t_i\big)=\\
  &=2-\big[ 2-(2-\varep_0)^2\big] \big(\sum_{i\in J_n}\ t_i\big)\le\\
  &\le 2-\big[ 2-(2-\varep_0)^2\big] \dfrac{1}{m}.
 \end{aligned}$$
\vk.2cm
Hence 

  $$2-\big[ 2-(2-\varep_0)^2\big] \dfrac{1}{m}>2-\dfrac{4}{n}$$

\no with fixed $\varep_0,\ m$ and for all $n$ satisfying (1), a contradiction.
\vk.2cm
Thus $\al(\{y_{q_{n-1}+1},\ y_{q_{n-1}+2},\ \cdots,\ y_{q_n}\}_{n=2}^\iy)=\sqrt{2}$. Since $d(A)=\sqrt{2}$ one concludes therefore $\al(A)=\sqrt{2}$.\hfill$\square$ 
\vk.3cm
For the second proof we need the following lemma which is a variation of 
[3, Lemma 4].
\vk.5cm
\textbf{Lemma\ 2.}\quad {\it Let $A$ be a non-empty bounded subset of Hilbert space $H$;\ \ $r$ and $c$ the Chebyshev radius of $A$ with respect to $H$ and Chebyshev center of $A$ in $H$, respectively. Then $c\in\ovl{co} A_\varep$ and $r=r(A_\varep)$ for every $\varep\in (0,r)$, where $A_\varep\col=A\setminus B(c,r-\varep)$.}
\vk.2cm
{\it Proof of Lemma 2.}\quad Assume contrariwise that $c$ is not the Chebyshev center of $A_\varep$ in $H$, then $r_1\col=r(A_\varep)<r$. Denoting by $c_1$ the Chebyshev center of $A_\varep$ in $H$ we choose $c'= \al c_1+(1-\al) c$ for some $\al\in (0,1)$ such that $0<\|c-c'\|<\varep$.
\vk.2cm
Take a point $x\in A$. If $x\in A_\varep$ then $\|x-c'\|\le \al\|x-c_1\|+(1-\al)\|x-c\|\le \al r_1+(1-\al) r<r$. In the other case $x\in A\setminus A_\varep$ we have $\|x-c'\|\le \|x-c\|+\|c-c'\|<r-\varep+\|c-c'\|<r$. 
\vk.2cm
So $A\subset B(c',r')$ with $r'<r$, a contradiction. The proof of the lemma is complete.\hfill$\square$
\vk.2cm
{\it The second proof of Theorem 1.}\quad In view of Lemma 2 by taking $\varep=\dfrac{1}{n}$ for every integer $n\ge 2 $ one has:\ \ $c\in\ovl{co} \Big(A\setminus B\big(c,1-\dfrac{1}{n}\big)\Big)$. Hence there exist $x_{p_{n-1}+1},\ x_{p_{n-1}+2},\ \cdots,\ x_{p_n}$ in $A\setminus B\big(c,1-\dfrac{1}{n}\big)$ and positive numbers $t_{p_{n-1}+1},\ t_{p_{n-1}+2},\ \cdots,\ t_{p_n}$ (with convention $p_1=0$) such that
  $$\sum_{p_{n-1}<i\le p_n}\ t_i=1,\qquad \Big\|\sum_{p_{n-1}<i\le p_n}\ t_i y_i-c\Big\|<\dfrac{1}{n}.$$

We show that $\al(\{x_{p_{n-1}+1},\ x_{p_{n-1}+2},\ \cdots,\ x_{p_n}\}_{n=2}^\iy)=\sqrt{2}$. Assume on the contrary that $\al(\{x_{p_{n-1}+1},\ x_{p_{n-1}+2},\ \cdots,\ x_{p_n}\}_{n=2}^\iy)<\sqrt{2}$. By choosing $\varep_0\in (0,\sqrt{2})$ satisfying inequality $\al(\{x_{p_{n-1}+1},\ x_{p_{n-1}+2},\ \cdots,\ x_{p_n}\}_{n=2}^\iy)\le\sqrt{2}-\varep_0$, there exist subsets $D_1,\ D_2,\ \cdots,\ D_m$ of $H$ with $d(D_i)\le\sqrt{2}-\varep_0$ for every $i=1, 2, \cdots, m$, such that $\{x_{p_{n-1}+1},\ x_{p_{n-1}+2},\ \cdots,\ x_{p_n}\}_{n=2}^\iy\subset\und{i=1}{\ov{m}{\bigcup}}\ D_i$. As in the first proof one can find among $D_1,\ D_2,\ \cdots,\ D_m$, a set, say $D_1$ with the property that there are infinitely many $n$ satisfying
  $$\sum_{i\in I_n}\ t_i\ge \dfrac{1}{m}\eqno{(2)}$$

\no where

  $$I_n\col=\{i\in [p_{n-1}+1,p_n]:\ \ x_i\in D_1\}.$$

Analogously for each $n$ satisfying (2) and fixed $j\in I_n$ we have

  $$\begin{aligned} & \sum_{p_{n-1}<i\le p_n}\ t_i \|x_i-x_j\|^2=\sum_{p_{n-1}<i\le p_n}\ t_i\|y_i-c+c-y_j\|^2=\\
  &{\quad} =\sum_{p_{n-1}<i\le p_n}\ t_i\big( \|y_i-c\|^2+\|y_j-c\|^2-2(y_i-c,y_j-c)\big)>\\
  &{\quad} >2\big(1-\dfrac{1}{n}\big)^2-2\big(\sum_{p_{n-1}<i\le p_n}\ t_i x_i-c,y_j-c\big)\ge\\
  &{\quad} \ge \big(1-\dfrac{1}{n}\big)^2-2\dfrac{1}{n}>2-\dfrac{6}{n}.
 \end{aligned}$$

Similarly one has also

  $$\sum_{p_{n-1}<i\le p_n}\ t_i \|x_i-x_j\|^2\le 2-\big[2-(\sqrt{2}-\varep_0)^2\big] \dfrac{1}{m}$$ 

Hence 

  $$2-\big[ 2-(2-\varep_0)^2\big] \dfrac{1}{m}>2-\dfrac{6}{n}$$

\no for all $n$ satisfying (2), a contradiction.
\vk.2cm
Thus $\al(\{x_{p_{n-1}+1},\ x_{p_{n-1}+2},\ \cdots,\ x_{p_n}\}_{n=2}^\iy)=\sqrt{2}$. This implies $\al(A)=\sqrt{2}$.\hfill$\square$
\vk.2cm
As an immediate consequence one obtains Gulevich's result mentioned in the Introduction.
\vk.5cm
\textbf{Corollary}\ (\cite{Gul}).\quad {\it Let $A$ be a relatively compact set in a Hilbert space with $d(A)>0$. Then $r(A)<\dfrac{1}{\sqrt{2}} d(A)$.}
\vk.5cm
\textbf{Remarks.}\quad 1.\ \ In \cite{Das} another proof of equality $J(H)=J_s(H)=\dfrac{1}{\sqrt{2}}$ was given by H. Steinlein. Essentially the heart of the proof is a relation between the Lifshitz characteristic and self-Jung constant of a Banach space $X$:\ \ $\varkappa_0(X)\le \big(J_s(X)\big)^{-1}$ which can be extended to the case of metric spaces with convex structure. We shall come back to this problem in a forthcoming paper.
\vk.2cm
2.\ \ By using Lemma 2 we see that $A_\varep$ is also an extremal set and $\al(A_\varep)=\sqrt{2}$ for every $\varep\in (0,1)$.
\vk.2cm
3.\ \ Although $\al(A_\varep)=\sqrt{2}$ for all $\varep\in (0,1)$, we may have $\ovl{co} A\cap S(c,1)=\emptyset$ ({\it cf.} \cite{Das}). The following question arises: what can be said about $\al\big(A\cap S(c,1)\big)$ when $A\cap S(c,1)\ne\emptyset$? The answer is:\ \ $\al\big(A\cap S(c,1)\big)$ can take arbitrary values in $[0,\sqrt{2}]$. Below we produce some examples.
\vk.2cm 
\textbf{Example 1.}\quad Let $\{e_n\}_{n=1}^\iy$ be an infinite orthonormal sequence in Hilbert space $H$. Set $A_1\col=\big\{\big(1-\dfrac{1}{n}\big) e_n\big\}_{n=1}^\iy$, and $A_2\col=\{x_1, x_2, \cdots, x_n, \cdots\}$ with

  \begin{align} & x_1\col=\dfrac{1}{\sqrt{2}} e_1+\dfrac{1}{\sqrt{2}} e_2;\quad x_2\col=\dfrac{1}{\sqrt{2}} e_1+\dfrac{1}{2} e_2+\dfrac{1}{2} e_3;\notag\\
  & x_n\col=\dfrac{1}{\sqrt{2}} e_1+\dfrac{1}{\sqrt{2^2}} e_2+ \cdots + \dfrac{1}{\sqrt{2^n}} e_n+\dfrac{1}{\sqrt{2^n}} e_{n+1};\ \cdots\notag
  \end{align}

It is easy to see that $r(A_1)=1,\ \ d(A_1)=\sqrt{2}$ and 0 is the Chebyshev center of $A_1$ in $H$. Furthermore $\|x_n\|=1$ for every $n$,\ $\|x_m-\big(1-\dfrac{1}{n}\big) e_n\|\le\sqrt{2},\ \ \for m, n$;\ \ $\|x_{n+p}-x_n\|^2=\dfrac{1}{2^n} \rg 0$ as $n\rg\iy$. Thus $\{x_n\}$ is a Cauchy sequence and one gets $\al(A_2)=0$.
\vk.2cm
Now setting $A\col=A_1\cup A_2$ we have $r(A)=1,\ \ d(A)=\sqrt{2}$ and 0 is also the Chebyshev center of $A$ in $H$. Obviously $A\cap S(0,1)=\ovl{co} A\cap S(0,1)=A_2$.
\vk.5cm
\textbf{Example 2.}\quad Let $\{e_n\}_{n=1}^\iy$ and $A_1$ be as in Example 1. For each $\ga\in (0,\sqrt{2}]$ putting $\be\col=\dfrac{\ga}{\sqrt{2}}\in (0,1]$ we choose $\la\in [0,1)$ satisfying $\la^2+\be^2=1$. Denote by $A_2\col=\{y_1, y_2, \cdots, y_n, \cdots\}$ with

  $$y_1\col=\la e_1+\be e_2;\ \ y_2\col=\la e_1+\be e_3;\ \cdots\ y_n\col=\la e_1+\be e_{n+1};\ \cdots$$

Obviously $\|y_n\|=1$ for every $n$;\ \ $\|y_n-y_m\|=\sqrt{2}\be=\ga,\ \ \for m\ne n$. Setting $A\col=A_1\cap A_2$ we obtain $r(A)=1,\ \ d(A)=\sqrt{2}$;\ \ 0 is the Chebyshev center of $A$ in $H$ and $A\cap S(0,1)=\ovl{co} A\cap S(0,1)=A_2$ with $\al(A_2)=\ga$.    
\vk.5cm
\begin{center}
\textbf {3.\ \ Proof of the Main Theorem}
\end{center}
\vk.2cm

From the first proof of Theorem 1 we derived a sequence $y_{q_{n-1}+1},\ y_{q_{n-1}+2},\ \cdots,\ y_{q_n}$ in $A_n\cap S(c_n,r_n)$\ (with $q_1=0$) and positive $t_{q_{n-1}+1},\ t_{q_{n-1}+2},\ \cdots,\ t_{q_n}$ for each integer $n\ge 2$ satisfying
  $$c_n=\sum_{q_{n-1}<i\le q_n}\ t_i y_i,\qquad \sum_{q_{n-1}<i\le q_n}\   t_i=1.$$
\vk.2cm
We claim that $\chi\big(\{y_{q_{n-1}+1},\ y_{q_{n-1}+2},\ \cdots,\ y_{q_n}\}_{n=2}^\iy\big)=1$. Assume that $A$ can be covered by a finite number of balls of radius $r$:\ \ $B_1,\ B_2,\ \cdots,\ B_m$. Then there exist a ball among $B_1,\ B_2,\ \cdots,\ B_m$, say $B_1$, such that there are infinitely many $n$ satisfying  

$$\sum_{i\in J_n}\ t_i\ge \dfrac{1}{m}\eqno{(3)}$$

\no where

  $$J_n\col=\{i\in [q_{n-1}+1,q_n]:\ \ y_i\in B_1\}.$$
\vk.2cm
As in the proof of Theorem 1 one has

  $$\sum_{q_{n-1}<i\le q_n}\ t_i \|y_i-y_j\|^2=2 r_n^2>2-\dfrac{4}{n}\eqno{(4)}$$

\no for every $j\in [q_{n-1}+1,q_n]$ fixed.
\vk.2cm
From (4) it follows that 

  $$\sum_{i\in I_{nj}}\ t_i<\dfrac{1}{\sqrt{n}},$$

\no where

  $$I_{nj}\col=\bigg\{i\in [q_{n-1}+1,q_n]:\ \ \|y_i-y_j\|^2<2-\dfrac{4}{\sqrt{n}}\bigg\},\ \ j\in [q_{n-1}+1,q_n],$$

\no and

  $$2(1-t_j)\ge \sum_{q_{n-1}<i\le q_n}\ t_i \|y_i-y_j\|^2=2 r_n^2>2-\dfrac{4}{n}.$$ 

This implies $t_j<\dfrac{2}{n}$ for every $j\in [q_{n-1}+1,q_n]$. Therefore if $n$ satisfies (3), then $|J_n|\dfrac{2}{n}>\dfrac{1}{m}$, or equivalently
 $|J_n|>\dfrac{n}{2m}$ (here $|J_n|$ denotes the cardinal of $J_n$).
\vk.2cm
For each $n$ satisfying (3) and $j\in J_n$ let us denote by 
  $$J_n(y_j)\col=\bigg\{i\in J_n:\ \ \|y_i-y_j\|^2\ge 2-\dfrac{4}{\sqrt{n}}\bigg\},$$

\no and
  
  $$\hat{J}_n(y_j)\col=\{y_i:\ \ i\in J_n(y_j)\}.$$  

Obviously from (5) one gets

  $$\sum_{i\in J_n\setminus J_n(y_j)}\ t_i<\dfrac{1}{\sqrt{n}}\eqno{(6)}$$

\no and

  $$\sum_{i\in J_n(y_j)}\ t_i>\dfrac{1}{m}-\dfrac{1}{\sqrt{n}}\eqno{(7)}$$

For each positive integer $p$ choose $n$ sufficiently large satisfying (3) and such that $\dfrac{p+1}{\sqrt{n}}\le\dfrac{1}{\sqrt{m}}$. We claim that for every choice of $i_1, i_2, \cdots, i_p\in J_n$ we have       $$\und{k=1}{\ov{p}{\bigcap}}\ J_n(y_{i_k})\ne\emptyset\eqno{(8)}$$ 

Indeed otherwise $\und{k=1}{\ov{p}{\bigcap}}\ J_n(y_{i_k})=\emptyset$ would imply that 

  $$J_n(y_{i_1})\subset J_n\setminus\big(\bigcap_{k=2}^p\ J_n(y_{i_k})\big)=\bigcup_{k=2}^p\ \big( J_n\setminus J_n(y_{i_k})\big).$$

\no Consequently by (6) and (7)

  $$\dfrac{1}{m}-\dfrac{1}{\sqrt{n}}<\sum_{\al\in J_n(y_{i_1})}\ t_\al\le
\sum_{k=2}^p\ \ \sum_{\al\in J_n\setminus J_n(y_{i_k})}\ t_\al<(p-1)\dfrac{1}{\sqrt{n}}.$$  

Thus $\dfrac{1}{m}<\dfrac{p}{\sqrt{n}}$. This would contradict to the choice of $n$ and $p$.
\vk.2cm
Next from (8) it follows that if $1\le k\le p$ and $i_1, i_2, \cdots, i_k\in J_n$, then $\und{\al=1}{\ov{k}{\bigcap}}\ \hat{J}_n(y_{i_\al})\ne\emptyset$. With $n$ and $p$ chosen as above let us fix $j\in J_n$. Setting $z_1\col=y_j$ we take consecutively $z_2\in \hat{J}_n(z_1);\ \ z_3\in\hat{J}_n(z_1)\cap\hat{J}_n(z_2);\ \ \cdots;\ \ z_{p+1}\in\und{i=1}{\ov{p}{\bigcap}}\ \hat{J}_n(z_i).$ 
\vk.2cm
Obviously $\|z_i-z_j\|^2\ge 2-\dfrac{4}{\sqrt{n}}$ for all $i\ne j$ in $\{1, 2, \cdots, p+1\}$. Now for given $\varep\in (0,\sqrt{2})$ choose $n$ as above and moreover sufficiently large so that $2-\dfrac{4}{\sqrt{n}}\ge (\sqrt{2}-\varep)^2$. One sees that $z_1, z_2, \cdots, z_{p+1}$ form a\ $p$-simplex $\De$ whose edges have length not less than $\sqrt{2}-\varep$.
\vk.2cm
We now prove that the radius $r$ of balls $B_1,\ B_2,\ \cdots,\ B_m$ is $\ge 1$. Let $c'$ and $r'$ denote respectively the Chebyshev center of $\De$ in $H$ and the Chebyshev radius of $\De$ with respect to $H$. From the proof of  classical Jung's theorem it follows that there exist non-negative $\al_1, \al_2, \cdots, \al_{p+1}$ with $\und{i=1}{\ov{p+1}{\sum}}\ \al_i=1$ and $c'=\und{i=1}{\ov{p+1}{\sum}}\ \al_i z_i$. Next for each $j\in\{1, 2, \cdots, p+1\}$ we have
\vk.2cm
  \begin{align} \big(2-\dfrac{4}{\sqrt{n}}\big)(1-\al_j)\ &  \le\sum_{i=1}^{p+1}\ \al_i \|z_i-z_j\|^2=\sum_{i=1}^{p+1}\ \al_i \|z_i-c'+c'-z_j\|^2=\notag\\
  &=\sum_{i=1}^{p+1}\ \al_i \big(\|z_i-c'\|^2+\|z_j-c'\|^2\big)
  -\bigg(\sum_{i=1}^{p+1}\ \al_i (z_i-c'),z_j-c'\bigg)\le\notag\\ 
  &\le 2(r')^2.\notag
  \end{align}
\vk.2cm
\no Thus

  $$\big(2-\dfrac{4}{\sqrt{n}}\big)\sum_{j=1}^{p+1}\ (1-\al_j)\le 2(p+1)(r')^2,$$

\no or equivalently

  $$r'\ge \sqrt{\dfrac{\big(2-4/\sqrt{n}\big)\ p}{2(p+1)}}\eqno{(9)}$$

The RHS of (9) tends to 1 as $p\rg\iy$. Obviously $r\ge r'$ since $\De\subset B_1$. This implies $r\ge 1$ as claimed. One concludes therefore $\chi(A)=1$.
\vk.2cm
Conversly if $d(A)=\sqrt{2}$, and for every $\varep\in (0,\sqrt{2})$ and every positive integer $p$\ \ $A$ contains a\ $p$-simplex $\De$ with its edges having length $\ge \sqrt{2}-\varep$, then we see immediately that $A$ is an extremal set.\hfill$\square$
\vk.5cm
\textbf{Acknowledgement.}\quad The second author would like to thank Professors Do Hong Tan and Nguyen Thanh Ha for encouragement and helpful discussions.
\vk.5cm


\begin{thebibliography}{1}
\bibitem{Ber}  Berdyshev, V. I., \emph {Connection between Jackson's inequality and a geometrical problem,} Math. Zametki, \textbf{3} (1968), 327--338 (Russian).
\bibitem{Byn} Bynum, W. L., \emph {Normal structure coefficients for Banach spaces,} Pacific J. Math. , \textbf{86} (1980), 427--436.
\bibitem{Das} Dane\v s, J., \emph {On the radius of a set in a Hilbert space,} Comment. Math. Univ. Carolina, \textbf{25}-2 (1984), 355--362.
\bibitem{Dgk} Danzer, L., Gr\"unbaum, B., Klee, V., \emph {Helly's theorem and its relatives,} Proc. Symp. Pure Math., Amer. Math. Soc., \textbf{7} (1963), 101--180. 
\bibitem{Gar} Garkavy, A. L., \emph {On the Chebyshev center and convex hull of a set,} Uspekhi Math. Nauk, \textbf{19} (1964), N$_-^0$ 6, 139--145. 
\bibitem{Gul} Gulevich, N. M., \emph {The radius of a compact set in a Hilbert space,} Zapiski Nauchnykh Seminarov LOMI im. V. A. Steklova, AN SSSR, vol. \textbf{164} (1988), 157--158 (Russian).
\bibitem{Jun} Jung, H. W. E., \emph {\"Uber die kleinste Kugel, die eine r\"aumliche Figur einschliesst,} J. Reine Angew. Math., \textbf{123} (1901),
241--257.  
\bibitem{Rou} Routledge, N. A., \emph {A result in Hilbert space,} Quart. J.  Math., Oxford, (2) \textbf{3} (1952), N$_-^0$ 9, 12--18.
\end{thebibliography}
\end{document}